\documentclass[12pt,leqno]{amsart} 
\newtheorem{thm}{Theorem}[section]
\newtheorem{defi}{Definition}[section]
\newtheorem{lem}{Lemma}[section]
\newtheorem{cor}{Corollary}[section]

\newtheorem{rem}{Remark}[section]
\newcommand{\be}{\begin{equation}}
\newcommand{\ee}{\end{equation}}
\newcommand{\bea}{\begin{eqnarray}}
\newcommand{\eea}{\end{eqnarray}}
\newcommand{\beb}{\begin{eqnarray*}}
\newcommand{\eeb}{\end{eqnarray*}}
\usepackage{amssymb,amsfonts,amsthm,setspace,indentfirst,mathrsfs}
\usepackage{minibox}
\numberwithin{equation}{section}
\begin{document}
%
\title[On warped product pseudosymmetric type manifolds]{\bf{On warped product manifolds satisfying some pseudosymmetric type conditions}}
\author[A. A. Shaikh and H. Kundu]{Absos Ali Shaikh and Haradhan Kundu}
\date{\today}
\address{\noindent\newline Department of Mathematics,\newline University of 
Burdwan, Golapbag,\newline Burdwan-713104,\newline West Bengal, India}
\email{aask2003@yahoo.co.in, aashaikh@math.buruniv.ac.in}
\email{kundu.haradhan@gmail.com}
\dedicatory{}
\begin{abstract}
The object of the present paper is to study the characterization of warped product manifolds satisfying some pseudosymmetric type conditions, especially, due to projective curvature tensor. For this purpose we consider a warped product manifold satisfying the pseudosymmetric type condition $R\cdot R = L_1 Q(g,R) + L_2 Q(S,R)$ and evaluate its characterization theorem. As special cases of $L_1$ and $L_2$ we find out the necessary and sufficient condition for a warped product manifold to satisfy various pseudosymmetric type, such as pseudosymmetry, Ricci generalized pseudosymmetry, semisymmetry due to projective curvature tensor ($P\cdot R = 0$), pseudosymmetry due to projective curvature tensor ($P\cdot R = L Q(g,R)$) etc. Finally we present some suitable examples of warped product manifolds satisfying such pseudosymmetric type conditions.
\end{abstract}
%
\maketitle
%

\section{Introduction}
%
\indent Let $\nabla$, $R$, $S$, $G$, $P$ and $\kappa$ be respectively the Levi-Civita connection, the Riemann-Christoffel curvature tensor, Ricci tensor, the Gaussian curvature tensor, the projective curvature tensor and the scalar curvature of an $n$-dimensional $(n\geq 3)$ connected smooth semi-Riemanian manifold $M$ equipped with the semi-Riemannian metric $g$. Symmetry is a very important geometric property of a space. Cartan \cite{Cart26} introduced the notion of symmetry (local and global) on a Riemannian manifold in terms of geodesic symmetries. A semi-Riemannian manifold is said to be locally symmetric \cite{Cart26} if its local geodesic symmetries at each point are all isometry. According to Cartan-Ambrose-Hicks theorem a locally symmetric manifold can be characterized by the curvature condition $\nabla R =0$, i.e., the curvature tensor is covariantly constant (For the meaning and definition of various  notations and symbols used here, we refer the reader to Section \ref{prel} of this paper).\\
\indent As a proper generalization of locally symmetric manifold, Cartan \cite{Cart46} introduced the notion of semisymmetric manifold. A semi-Riemannian manifold is said to be semisymmetric \cite{Cart46} (see also \cite{Szab82}, \cite{Szab84}, \cite{Szab85}) if $R\cdot R = 0$, where the first $R$ stands for the curvature operator acting as a derivation on the second $R$. It may be noted that a semi-Riemannian manifold is semisymmetric if and only if its sectional curvature function $k(p,\pi)$ is invariant up to second order, under parallel transport of any plane $\pi$ at any point $p$ of $M$ around any infinitesimal coordinate parallelogram centered at $p$ (see \cite{DHV08}, \cite{HV06}, \cite{HV07}).\\
\indent During the study of totally umbilical submanifolds of semisymmetric manifolds as well as during the consideration of geodesic mappings on semisymmetric manifolds, Adamow and Deszcz \cite{AD83} (see \cite{Desz92} and also references therein) introduced the notion of pseudosymmetric manifolds as a proper generalization of semisymmetric manifolds. A semi-Riemannian manifold is said to be pseudosymmetric if 
$$
\mbox{$R\cdot R$ and $Q(g,R)$ are linearly dependent.}
$$
Let $M$ be not of constant curvature and $U$ denotes the set $\{x \in M: Q(g, R) \neq 0\}$. Then at a point $p \in U$, a plane $\pi_1(\vec{v_p}, \vec{w_p}) \subset T_p M$ is said to be curvature dependent with respect to another plane $\pi_2(\vec{x_p}, \vec{y_p}) \subset T_pM$ if $Q(g, R)(\vec{v_p}, \vec{w_p}, \vec{v_p}, \vec{w_p}, \vec{x_p}, \vec{y_p}) \neq 0$. Now if $\pi_1$ is curvature dependent with respect to $\pi_2$, then the scalar
$$L(p,\pi_1,\pi_2) = \frac{R\cdot R(\vec{v_p}, \vec{w_p}, \vec{v_p}, \vec{w_p}, \vec{x_p}, \vec{y_p})}{Q(g,R)(\vec{v_p}, \vec{w_p}, \vec{v_p}, \vec{w_p}, \vec{x_p}, \vec{y_p})}$$
is called the double sectional curvature or Deszcz sectional curvature (\cite{DHV08}, \cite{HV06}, \cite{HV07}, \cite{HV09}) of the plane $\pi_1$ with respect to $\pi_2$ at $p$. In terms of Deszcz sectional curvature, a semi-Riemannian manifold is pseudosymmetric if at each point $p\in U$, $L(p,\pi_1,\pi_2)$ is independent of the planes $\pi_1$ and $\pi_2$ (see \cite{DHV08}, \cite{HV06}, \cite{HV07}, \cite{HV09}).\\
\indent Replacing $R$ by other curvature tensors and $g$ by other symmetric $(0,2)$-tensor in the defining condition of semisymmetric manifold and pseudosymmetric manifold one can get various curvature restricted geometric structures, which are simply called as semisymmetric type and pseudosymmetric type manifolds (\cite{DGHZ15}, \cite{DGJZ16}, \cite{DHJKS14}, \cite{SK14}, \cite{SKsrs}, \cite{SRK15}). For the geometric meaning of Weyl semisymmetric and Weyl pseudosymmetric spaces we refer the reader to see \cite{JHPV09}. One of the important semisymmetric type (resp., pseudosymmetric type) manifold is semisymmetric (resp., pseudosymmetric) manifold due to projective curvature tensor. A semi-Riemannian manifold is said to be semisymmetric (resp., pseudosymmetric) due to projective curvature tensor if 
\begin{center}
$P\cdot R = 0$ (resp., $P\cdot R$ and $Q(g,R)$ are linearly dependent).
\end{center}
Another important pseudosymmetric type manifold is Ricci generalized pseudosymmetric manifold. A semi-Riemannian manifold is said to be Ricci generalized pseudosymmetric (\cite{DD91}, \cite{DD91a}) if 
\begin{center}
$R\cdot R$ and $Q(S,R)$ are linearly dependent.
\end{center}
We refer the reader to see \cite{SKppsn} for details about various curvature restricted geometric structures due to projective curvature tensor.\\
\indent Again the notion of warped product manifold (\cite{BO69}, \cite{Kr57}) is a generalization of product manifold and this notion is important due to its applications in general theory of relativity and cosmology. Various spacetimes are warped product, e.g., Robertson-Walker spacetimes, asymptotically flat spacetimes, Schwarzschild spacetimes, Kruskal space-times, Reissner-Nordstr\"{o}m spacetimes etc. The main purpose of this present paper is to study the characterization of a warped product semi-Riemannian manifold realizing some pseudosymmetric type curvature conditions, especially, due to the projective curvature tensor. For this purpose we consider the pseudosymmetric type condition $R\cdot R = L_1 Q(g,R) + L_2 Q(S,R)$, and evaluate the characterization of a warped product manifold satisfying such curvature condition. As a special case we get the characterization of a warped product manifold which is (i) semisymmetric, (ii) pseudosymmetric, (iii) Ricci generalized pseudosymmetric, (iv) special Ricci generalized pseudosymmetric, (v) semisymmetric due to projective curvature tensor and (vi) pseudosymmetric due to projective curvature tensor etc. It is shown that if a warped product manifold $M = \overline M \times_f \widetilde M$ satisfies $R\cdot R = L_1 Q(g,R) + L_2 Q(S,R)$ such that $L_2$ is nowhere zero, then at either the base $\overline M$ is flat or the fiber $\widetilde M$ is Einstein. Consequently for a warped product pseudosymmetric manifold due to projective curvature tensor or special Ricci generalized pseudosymmetric manifold either the base is flat or the fiber is Einstein.\\
\indent The paper is organized as follows. After discussing various notations as preliminaries in Section 2, we define various pseudosymmetric type curvature restricted geometric structures in Section 3. Section 4 is devoted to the study of warped product manifold and we state the curvature relation of a warped product manifold with its base and fiber. In Section 5 we discuss about the characterization theorems of various pseudosymmetric type warped product manifolds. Finally to support our results we present some suitable examples of warped product manifolds in the last section. It is interesting to mention that notion of pseudosymmetric manifold arose during the study of totally umbilical hypersurface of a semisymmetric manifold, and in Example 1 we present a pseudosymmetric totally umbilical hypersurface of a semisymmetric manifold.
%
\section{Preliminaries}\label{prel}
Let $M$ be a connected $n$-dimensional smooth manifold equipped with the semi-Riemannian metric $g$. Let us consider the following notations related to $(M,g)$:\\
\indent $C^{\infty}(M) =$ the algebra of all smooth functions on $M$, \\
\indent $\mathcal T^r_k(M) =$ the space of all smooth tensor fields of type $(r,k)$ on $M$ and\\
\indent $\chi(M) = \mathcal T^1_0(M) =$ the Lie algebra of all smooth vector fields on $M$.\\
\indent The Kulkarni-Nomizu product (\cite{DGHS11}, \cite{Glog02}, \cite{SRK16}) $A\wedge E \in \mathcal T^0_4(M)$ of $A$ and $E \in \mathcal T^0_2(M)$, is given by
\beb
(A \wedge E)(X_1,X_2,X_3,X_4)&=&A(X_1,X_4)E(X_2,X_3) + A(X_2,X_3)E(X_1,X_4)\\\nonumber
&-&A(X_1,X_3)E(X_2,X_4) - A(X_2,X_4)E(X_1,X_3),
\eeb
where $X_1,X_2,X_3,X_4 \in \chi(M)$. Throughout the paper we consider $X,Y,X_1,X_2,\cdots \in \chi(M)$.\\
\indent Now for $D\in \mathcal T^0_4(M)$, $A\in \mathcal T^0_2(M)$ and $X,Y\in \chi(M)$, we get two endomorphisms $\mathscr D(X,Y)$ (called the associated curvature operator of $D$) and $X\wedge_A Y$ defined by
$$\mathscr{D}(X,Y)(X_1)=\mathcal D(X,Y)X_1 \ \mbox{and}$$
$$(X\wedge_A Y)X_1 = A(Y,X_1)X - A(X,X_1)Y,$$
where $\mathcal D \in \mathcal T^1_3(M)$ such that $g(\mathcal{D}(X,Y)X_1,X_2)=D(X,Y,X_1, X_2)$, called the associated $(1,3)$ tensor of $D$.\\
\indent A tensor $D\in \mathcal T^0_4(M)$ is said to be a generalized curvature tensor (\cite{DGHS11}, \cite{SDHJK15}, \cite{SK14}) if
$$D(X_1,X_2,X_3,X_4)+D(X_2,X_3,X_1,X_4)+D(X_3,X_1,X_2,X_4)=0,$$
$$D(X_1,X_2,X_3,X_4)+D(X_2,X_1,X_3,X_4)=0 \ \ \mbox{and}$$
$$D(X_1,X_2,X_3,X_4)=D(X_3,X_4,X_1,X_2).$$
The Gaussian curvature tensor $G$, Weyl conformal curvature tensor $C$, concircular curvature tensor $W$ and conharmonic curvature tensor $K$are all generalized curvature tensors and respectively given by
$$G = \frac{1}{2} g \wedge g,$$
$$C = R -\frac{1}{n-2} g \wedge S+\frac{\kappa}{2(n-1)(n-2)}g \wedge g,$$ 
$$W = R - \frac{\kappa}{2 n(n-1)}g \wedge g \ \ \mbox{and}$$ 
$$K = R - \frac{1}{n-2} g \wedge S.$$
The projective curvature tensor $P$ of type $(0, 4)$ , given by 
$$P(X_1,X_2,X_3, X_4) = R(X_1,X_2,X_3, X_4) -\frac{1}{n-2} \left[S(X_2,X_3)g(X_1,X_4)-S(X_1,X_3)g(X_2,X_4)\right],$$ 
is not a generalized curvature tensor.\\
\indent Again an endomorphism $\mathscr L$ can be operate on a $(0,k)$-tensor $H$ and obtain $\mathscr{L} H$ as follows:
\beb
(\mathscr{L} H)(X_1,X_2,\cdots,X_k) &=& -H(\mathscr{L}X_1,X_2,\cdots,X_k) - \cdots - H(X_1,X_2,\cdots,\mathscr{L}X_k).
\eeb
In particular, for $\mathscr L = \mathscr D(X,Y)$ and $X\wedge_A Y$ we get two $(0,k+2)$ tensors $D\cdot H$ and $Q(A,H)$ defined as (\cite{SKgrt}, \cite{SRK15}, \cite{DGJZ16} and also references therein)
\beb
&&D\cdot H(X_1,X_2, \ldots ,X_k,X,Y) = (\mathscr{D}(X,Y)\cdot H)(X_1,X_2, \ldots ,X_k)\\
&&= -H(\mathcal{D}(X,Y)X_1,X_2, \ldots ,X_k) - \cdots - H(X_1,X_2, \ldots ,\mathcal{D}(X,Y)X_k),
\eeb
\beb
&&Q(A,H)(X_1,X_2, \cdots ,X_k,X,Y) = ((X \wedge_A Y)\cdot H)(X_1,X_2, \ldots ,X_k)\\
&&= A(X, X_1) H(Y,X_2,\cdots,X_k) + \cdots + A(X, X_k) H(X_1,X_2,\cdots,Y)\\
&&- A(Y, X_1) H(X,X_2,\cdots,X_k) - \cdots - A(Y, X_k) H(X_1,X_2,\cdots,X).
\eeb
\begin{lem}\label{lem1}
Let $M$ be a connected warped product manifold with base $\overline M$ and fiber $\widetilde M$. If $f_1 \in C^{\infty} (\overline M)$ and $f_2 \in C^{\infty} (\widetilde M)$ satisfies $f_1~f_2\equiv 0$, then either $f_1\equiv 0$ or $f_2\equiv 0$.
\end{lem}
\begin{lem}\label{lem2}\cite{DD91} 
Let $A,E\in\mathcal T^0_2(M)$ be symmetric. Then $Q(A,E)=0$ if and only if $A$ and $E$ are linearly dependent \cite{DD91}. Moreover $Q(A,E)=0$ if and only if $M$ is Einstein, i.e., $S = \frac{\kappa}{n} g$.
\end{lem}
\section{Pseudosymmetric type curvature restricted geometric structures}
\begin{defi} (\cite{SK14}, \cite{Szab82}, \cite{Szab84}, \cite{Szab85}) 
For $H \in \mathcal T^0_k(M)$ and $D \in \mathcal T^0_4(M)$, a semi-Riemannian manifold $M$ is said to be $H$-semisymmetric type if $D\cdot H = 0$.
\end{defi}
In particular, a semi-Riemannian manifold with the semisymmetric type conditions 
$R\cdot R = 0$, 
$R\cdot S = 0$, 
$R\cdot P = 0$, 
$P\cdot R = 0$ and 
$P\cdot S = 0$
is respectively called semisymmetric, Ricci semisymmetric, projective semisymmetric, semisymmetric due to projective curvature tensor and Ricci semisymmetric due to projective curvature tensor respectively.
\begin{defi}$($\cite{AD83}, \cite{Desz92}, \cite{DGHS11}, \cite{SK14}$)$ 
For $H \in \mathcal T^0_k(M)$ and $D_i \in \mathcal T^0_4(M)$, $i=1,2,\cdots r, r\ge 2$, a semi-Riemannian manifold is said to be $H$-pseudosymmetric type if $D_i\cdot H$ are linearly dependent, i.e., $\sum\limits_{i=1}^r c_i (D_i\cdot H) = 0$ for some $c_i\in C^{\infty}(M)$, called the associated scalars. Moreover if the associated scalars are constant then the manifold is called pseudosymmetric type manifold of constant type.
\end{defi}
In particular, a semi-Riemannian manifold satisfying 
$$R\cdot C = L_R Q(g,C) \ \ \mbox{on} \ \ \left\{x\in M: C_x \ne 0\right\},$$ 
$$R\cdot S = L_S Q(g,S) \ \ \mbox{on} \ \ \left\{x\in M: \left(S-\frac{\kappa}{n}g\right)_x \ne 0\right\},$$
$$R\cdot P = L_P Q(g,P) \ \ \mbox{on} \ \ \left\{x\in M: Q(g,P)_x \ne 0\right\},$$ 
$$P\cdot R = L_1 Q(g,R) \ \ \mbox{on} \ \ \left\{x\in M: \left(R-\frac{\kappa}{n(n-1)}G\right)_x \ne 0\right\} \ \mbox{and}$$ 
$$P\cdot S = L_2 Q(g,S) \ \ \mbox{on} \ \ \left\{x\in M: \left(S-\frac{\kappa}{n}g\right)_x \ne 0\right\}$$
is called conformally pseudosymmetric, Ricci pseudosymmetric, projective pseudosymmetric, pseudosymmetric due to projective curvature tensor and Ricci pseudosymmetric due to projective curvature tensor respectively, where $L_R, L_S,L_P,L_1,L_2$ are the associated scalars.\\
\indent In this paper we are mainly interested on the pseudosymmetric type condition $R\cdot R = L_1 Q(g,R) + L_2 Q(S,R)$. For particular values of $L_1$ and $L_2$, we get various pseudosymmetric type conditions as follows:
\begin{center}
\begin{tabular}{||c|c|c||}\hline\hline
Values of $L_1$ \& $L_2$ & Structure name & Curvature condition\\\hline

$L_1 = L_2 = 0$ & Semisymmetric manifold & $R\cdot R = 0$\\\hline

$L_2 = 0$ & Pseudosymmetric manifold & $R\cdot R = L_1 Q(g,R)$\\\hline

$L_1 = 0$ & $\begin{array}{c} \mbox{Ricci generalized}\\ \text{Pseudosymmetric manifold} \end{array}$ & $R\cdot R = L_2 Q(S,R)$\\\hline

$L_1 = 0, L_2 = 1$ & $\begin{array}{c} \mbox{Special Ricci generalized}\\ \text{Pseudosymmetric manifold} \end{array}$ & $R\cdot R = Q(S,R)$\\\hline

$L_1 = \frac{\kappa}{n(n-1)}, L_2 = 0$ & $\begin{array}{c} \mbox{Semisymmetric due to}\\ \text{concircular curvature tensor} \end{array}$ & $W\cdot R = 0$\\\hline

$L_1 = \frac{\kappa}{n(n-1)}$ & $\begin{array}{c} \mbox{Ricci generalized Pseudosymmetric}\\ \text{due to concircular curvature tensor} \end{array}$ & $W\cdot R = L_2 Q(S,R)$\\\hline

$L_1 = 0, L_2 = \frac{1}{n-2}$ & $\begin{array}{c} \mbox{Semisymmetric due to}\\ \text{projective curvature tensor} \end{array}$ & $P\cdot R = 0$\\\hline

$L_2 = \frac{1}{n-2}$ & $\begin{array}{c} \mbox{Pseudosymmetric due to}\\ \text{projective curvature tensor} \end{array}$ & $P\cdot R = L_1 Q(g,R)$\\\hline\hline
\end{tabular}
\end{center}
\section{Warped product manifolds}
Warped product is an important notion in semi-Riemannian geometry and it has a great applications in general theory of relativity and cosmology. To impose a semi-Riemannian structure on a product smooth manifold, the notion of warped product metric arose as a generalization of Riemannian product metric. This notion was independently introduced by Kru$\breve{\mbox{c}}$kovi$\breve{\mbox{c}}$ \cite{Kr57} (as semi-decomposable metric) as well as Bishop and O'Neill \cite{BO69}. The most important example of non-Riemann product but a warped product is surface of revolution. Many well-known spacetimes, e.g., Schwarzschild, Kottler, Reissner-Nordstr\"om, Reissner-Nordstr\"om-de Sitter, Vaidya, Robertson-Walker spacetimes are all warped products.\\
\indent Let $(\overline M, \overline g)$ and $(\widetilde M, \widetilde g)$ be two semi-Riemannian manifolds of dimension $p$ and $(n-p)$ respectively ($1\le p \le n-1$) and $M = \overline M \times \widetilde M$. The warped product metric $g$ on $M$ is given by
$$g=\pi^*(\overline g) + (f\circ\pi) \sigma^* (\widetilde g),$$
where $f$ is a positive smooth  function on $\overline M$ and $\pi:M\rightarrow\overline M$ (resp., $\sigma:M\rightarrow\widetilde M$) is the natural projection on $\overline M$ (resp., $\widetilde M$). The manifold $\overline M$ is called the base, $\widetilde M$ is called the fiber and $f$ is called the warping function of $M$. If $f=1$, then the warped product reduces to the Riemann product. If we consider a product chart $\left(U \times V; x^1,x^2, ..., x^p,x^{p+1}=y^1,x^{p+2}=y^2, ...,x^{n}=y^{n-p}\right)$ on $M$, then in terms of local coordinates, $g$ can be expressed as
\begin{eqnarray}\label{eq3.1}
g_{ij}=\left\{\begin{array}{lll}
&\overline g_{ij}&\ \ \ \ \mbox{for} \ i = a \ \mbox{and} \ j = b,\\
&f \widetilde g_{ij}&\ \ \ \ \mbox{for $i = \alpha$ and $j = \beta$,}\\
&0&\ \ \ \ \mbox{otherwise,}\\
\end{array}\right.
\end{eqnarray}
where $a,b \in \{1,2,...,p\}$ and $\alpha, \beta \in \{p+1,p+2,...,n\}$. We note that throughout the paper we consider $a,b,c,d,e,s,t\in \{1,2, ..., p\}$; $\alpha,\beta,\gamma,\delta,\epsilon,\mu,\eta\in \{p+1,p+2,...,n\}$ and $i,j,k,l,q,u,v\in \{1,2,...,n\}$. Moreover, when {\scriptsize $\Box$} is a quantity formed with respect to $g$, we denote by {\scriptsize $\overline \Box$} and {\scriptsize $\widetilde \Box$}, the similar quantities formed with respect to $\overline g$ and $\widetilde g$ respectively.\\
\indent By a straightforward calculation we can evaluate the components of various necessary tensors of a warped product manifold in terms of the base and fiber components. The non-zero local components $R_{hijk}$ of the Riemann-Christoffel curvature tensor $R$, $S_{jk}$ of the Ricci tensor $S$ and the scalar curvature $\kappa$ of $M$ are respectively given by
\beb
R_{abcd} = \overline{R}_{abcd},\,\,\,\, R_{a\alpha b\beta}=f T_{ab}\widetilde{g}_{\alpha \beta},\,\,\,R_{\alpha \beta \gamma \delta} = f\widetilde{R}_{\alpha \beta \gamma \delta} - f^2 \Delta \widetilde{G}_{\alpha \beta \gamma \delta}, 
\eeb
\beb
S_{ab}=\overline{S}_{ab}-(n-p)T_{ab},\,\,\,\, S_{\alpha \beta}=\widetilde{ S}_{\alpha \beta} + \Omega \widetilde{g}_{\alpha \beta}
\eeb
\beb
\mbox{and } \kappa=\overline{\kappa}+\frac{\widetilde{\kappa}}{f}-(n-p)[(n-p-1)\Delta - 2 \; tr(T)], 
\eeb
where $G_{ijkl} = g_{il}g_{jk}-g_{ik}g_{jl}$ are the components of Gaussian curvature and
$$T_{ab} = \frac{1}{2f}(f_{a,b} - \frac{1}{2f}f_a f_b), \ \ \ \ tr(T) = g^{ab}T_{ab},$$
$$\Delta = \frac{1}{4f^2}g^{ab}f_a f_b, \ \ \ \ \Omega = -f((n-p-1)\Delta +tr(T)), \ \ \ f_{a}=\partial_{a} f=\frac{\partial f}{\partial x^{a}}.$$
The non-identically zero local components of $R\cdot R$, $Q(g,R)$ and $Q(S,R)$ of $M$ are given by
\bea\label{wp-r.r}
(R\cdot R)_{abcdst} &=& (\overline R\cdot \overline R)_{abcdst},\\\nonumber
(R\cdot R)_{a\alpha b\beta st} &=& f \widetilde g_{\alpha\beta}(\overline R\cdot T)_{abst},\\\nonumber
(R\cdot R)_{abc\alpha s\eta} &=& f \widetilde g_{\alpha\eta}\left(T_{as}T_{bc}-T_{ac}T_{bs}+T^t_{s}\overline R_{abct}\right),\\\nonumber
(R\cdot R)_{a\alpha\beta\gamma s\eta} &=& -f T_{as} [\widetilde R_{\eta\alpha\beta\gamma}-f\Delta \widetilde G_{\eta\alpha\beta\gamma}] - f^2 T^2_{as} \widetilde G_{\eta\alpha\beta\gamma},\\\nonumber
(R\cdot R)_{\alpha\beta\gamma\delta\mu\eta} &=& f [(\widetilde R\cdot \widetilde R)_{\alpha\beta\gamma\delta\mu\eta}-f \Delta Q(\widetilde g,\widetilde R)_{\alpha\beta\gamma\delta\mu\eta}];
\eea
\bea\label{wp-qgr}
Q(g, R)_{abcdst} &=& Q(\overline g, \overline R)_{abcdst},\\\nonumber
Q(g, R)_{a\alpha b\beta st} &=& f \widetilde g_{\alpha\beta} Q(\overline g, T)_{abst},\\\nonumber
Q(g, R)_{abc\alpha s\eta} &=& f \widetilde g_{\alpha\eta}\left(\overline R_{abct} - g_{bs}T_{ac}+ g_{as}T_{bc}\right),\\\nonumber
Q(g, R)_{a\alpha\beta\gamma s\eta} &=& -f [\overline g_{as} \widetilde R_{\eta\alpha\beta\gamma} + f(T_{as}-\Delta \overline g_{as})\widetilde G_{\eta\alpha\beta\gamma}],\\\nonumber
Q(g, R)_{\alpha\beta\gamma\delta\mu\eta} &=& f^2 Q(\widetilde g,\widetilde R)_{\alpha\beta\gamma\delta\mu\eta};
\eea
\bea\label{wp-qsr}
Q(S, R)_{abcdst} &=& Q(\overline S, \overline R)_{abcdst} - (n-p)Q(\overline S, \overline R)_{abcdst},\\\nonumber
Q(S, R)_{a\alpha b\beta st} &=& f \widetilde g_{\alpha\beta} Q(\overline S, T)_{abst},\\\nonumber
Q(S, R)_{abc\alpha s\eta} &=& \overline R_{abcs} [\widetilde S_{\alpha\beta} + \Omega \widetilde g_{\alpha\beta}]\\\nonumber
							&&+ f \widetilde g_{\alpha\eta} [T_{bc}(\overline S_{as}-(n-p) T_{as}) - T_{ac}(\overline S_{bs}-(n-p) T_{bs})],\\\nonumber
Q(S, R)_{a\alpha\beta\gamma s\eta} &=& - f (\overline S_{as}-(n-p) T_{as})(\widetilde R_{\eta\alpha\beta\gamma}- f\Delta \widetilde G_{\eta\alpha\beta\gamma})\\\nonumber
																	&&- f T_{as} \left[\widetilde g_{\alpha\beta}(\widetilde S_{\gamma\eta} + \Omega \widetilde g_{\gamma\eta}) -\widetilde g_{\alpha\gamma}(\widetilde S_{\beta\eta} + \Omega \widetilde g_{\beta\eta})\right],\\\nonumber
Q(S, R)_{a\alpha b\beta\mu\eta} &=& - f T_{ab} Q(\widetilde g, \widetilde S)_{\alpha\beta\mu\eta},\\\nonumber
Q(S, R)_{\alpha\beta\gamma\delta\mu\eta} &=& f [Q(\widetilde S,\widetilde R)_{\alpha\beta\gamma\delta\mu\eta} - \Delta Q(\widetilde S,\widetilde G)_{\alpha\beta\gamma\delta\mu\eta} + f \Omega Q(\widetilde g,\widetilde R)_{\alpha\beta\gamma\delta\mu\eta}].
\eea
We refer the readers to see \cite{DD91}, \cite{SK12}, \cite{SKgrtw}, \cite{SKA16} and also references therein for detail information about warped product components of various tensors on $M$.
\section{Main results}
\begin{thm}\label{main}
Let $M^n = \overline M^p \times_f \widetilde M^{n-p}$ be a warped product manifold. Then $M$ satisfies the pseudosymmetric type condition
\be\label{eq4.1}
R\cdot R = L_1 Q(g,R) + L_2 Q(S,R)
\ee
if and only if the following conditions hold simultaneously:\\
\indent $(I)$ $\overline R\cdot \overline R = L_1 Q(\overline g,\overline R) + L_2 Q(\overline S,\overline R) - L_2 (n-p) Q(T,\overline R),$\\
\indent $(II)$ $f \widetilde g_{\alpha\beta} \left(T_{as}T_{bc}-T_{ac}T_{bs}+T^t_{s}\overline R_{abct}\right) = 
L_1 f  \widetilde g_{\alpha\beta} \left(\overline R_{abcs}-T_{ac}\overline g_{bs}+T_{bc}\overline g_{as}\right)
+L_2 \overline R_{abcs}\left(\widetilde S_{\alpha\beta} + \Omega \widetilde g_{\alpha\beta}\right)\\
\indent\hspace{2.6in} +L_2 f \widetilde g_{\alpha\beta}  \left[T_{bc} \left(\overline S_{as}-(n-p)T_{as}\right)-T_{ac} \left(\overline S_{bs}-(n-p)T_{bs}\right)\right]$,\\
\indent $(III)$ $T_{as}\left(\widetilde R_{\eta\alpha\beta\gamma}- f \Delta\widetilde G_{\eta\alpha\beta\gamma}\right) + f T^2_{as}\widetilde G_{\eta\alpha\beta\gamma}=
 L_1 \left[\overline g_{as}\widetilde R_{\eta\alpha\beta\gamma}+ f (T_{as}-\Delta\overline g_{as})\widetilde G_{\eta\alpha\beta\gamma}\right]\\
\indent\hspace{3.3in} + L_2 (\overline S_{as}-(n-p) T_{as})(\widetilde R_{\eta\alpha\beta\gamma}- f\Delta \widetilde G_{\eta\alpha\beta\gamma})\\
\indent\hspace{3.3in} + L_2 T_{as} \left[\widetilde g_{\alpha\beta}(\widetilde S_{\gamma\eta} + \Omega \widetilde g_{\gamma\eta}) -\widetilde g_{\alpha\gamma}(\widetilde S_{\beta\eta} + \Omega \widetilde g_{\beta\eta})\right]$,\\
\indent	$(IV)$ $L_2 T_{ab} Q(\widetilde g, \widetilde S)_{\alpha\beta\gamma\delta} =0$ and\\
\indent	$(V)$ $\widetilde R\cdot \widetilde R = f \left(\Delta + L_1 + \Omega L_2\right) Q(\widetilde g,\widetilde R) + L_2 Q(\widetilde S,\widetilde R) - L_2 \Delta Q(\widetilde S,\widetilde G)$.
\end{thm}
\noindent \textbf{Proof:} In terms of local coordinates, \eqref{eq4.1} can be written as 
\be\label{eq4.1l}
(R\cdot R)_{ijkluv} = L_1 Q(g,R)_{ijkluv} + L_2 Q(S,R)_{ijkluv}.
\ee
Now from \eqref{wp-r.r}, \eqref{wp-qgr} and \eqref{wp-qsr} we see that the non-zero possibilities of \eqref{eq4.1l} are\\
(i) $i=a,j=b,k=c,l=d,u=s,v=t$,\\
(ii) $i=a,j=\alpha,k=b,l=\beta,u=s,v=t$,\\
(iii) $i=a,j=b,k=c,l=\alpha,u=s,v=\eta$,\\
(iv) $i=a,j=\alpha,k=\beta,l=\gamma,u=s,v=\eta$,\\
(v) $i=a,j=\alpha,k=b,l=\beta,u=\mu,v=\eta$ and\\
(vi) $i=\alpha,j=\beta,k=\gamma,l=\delta,u=\mu,v=\eta$.\\
So it is obvious that to prove the theorem we have only to show that the conditions (I)-(V) are the simplified form of the possibilities (i)-(vi). Putting (i), (iii), (iv), (v) and (vi) in \eqref{eq4.1l} and simplifying we get (I) to (V) respectively. Again putting (ii) in \eqref{eq4.1l}, we get $\overline R\cdot T = L_1 Q(\overline g,T) + L_2 Q(\overline S,T)$, which obviously follows from (II). This completes the proof.
\begin{cor}
If a warped product manifold $M^n = \overline M^p \times_f \widetilde M^{n-p}$ satisfies the pseudosymmetric type condition
$R\cdot R = L_1 Q(g,R) + L_2 Q(S,R)$, then $\overline R\cdot T = L_1 Q(\overline g, T) + L_2 Q(\overline S, T)$.
\end{cor}
\noindent \textbf{Proof:} From condition (II) of Theorem \ref{main} we get the result easily.
\begin{thm}
If $M^n = \overline M^p \times_f \widetilde M^{n-p}$ is a warped product manifold satisfying \eqref{eq4.1}, then 
$$M= \{x\in M : \overline R|_{\pi(x)} =0\} \cup \{x\in M : (T-L_1 \overline g)|_{\pi(x)} =0\} \cup \{x\in M : (\widetilde S-\mbox{$\frac{\widetilde \kappa}{n-p}$}\widetilde g)|_{\sigma(x)} =0\}.$$
\end{thm}
\noindent \textbf{Proof:} Since $M$ satisfies \eqref{eq4.1}, then from the condition (IV) of Theorem \ref{main}, at every point $x\in M$ we have the following three cases:\\
Case 1: $L_2 = 0$. Therefore from (III) of Theorem \ref{main}, we get
$$(T-L_1 g)\widetilde R = f\left[L_1(T - \Delta \overline g)+ \Delta T - T^2\right] \widetilde G$$
$\Rightarrow$ either $T-L_1 \overline g = 0$ or $\widetilde R$ is a scalar multiple of $\widetilde G$ and thus $\widetilde S = \frac{\widetilde \kappa}{n-p} \widetilde g$ at $x$.\\
Case 2: $T =0$. In this case putting the value of $T$ in (II) of Theorem \ref{main}, we get
$$\overline R (L_1 f \widetilde g + L_2 \widetilde S + L_2 \Omega \widetilde g) = 0$$
$$\Rightarrow \mbox{ either } \overline R = 0  \mbox{ or } L_2 \widetilde S + (L_1 f+  L_2 \Omega) \widetilde g = 0 \mbox{ at } x$$
$$\Rightarrow \mbox{ either }  \overline R = 0 \mbox{ or }  L_2 = 0 \mbox{ or }  \widetilde S = \frac{\widetilde \kappa}{n-p} \widetilde g \mbox{ at }  x.$$
Case 3: $Q(\widetilde g, \widetilde S) = 0$. Then from Lemma 2.1 of \cite{DD91}, $\widetilde S = \frac{\widetilde \kappa}{n-p} ~\widetilde g$ at $x$.\\
Now $\overline R = 0$ at $x$ means $\overline R|_{\pi(x)} = 0$, $T-L_1 \overline g = 0$ at $x$ means $(T-L_1 \overline g)|_{\pi(x)} = 0$ and $\widetilde S-\frac{\widetilde \kappa}{n-p}\widetilde g = 0$ at $x$ means $(\widetilde S-\frac{\widetilde \kappa}{n-p}\widetilde g)|_{\sigma(x)} = 0$. Now combining the resulting condition of above cases we get our assertion.
\begin{thm}\label{thm5.3}
If a warped product manifold $M^n = \overline M^p \times_f \widetilde M^{n-p}$ satisfies the pseudosymmetric type condition \eqref{eq4.1} and $L_2$ is nowhere zero, then\\
either (i) the base $\overline M$ is flat\\
or (ii) the fiber $\widetilde M$ is Einstein.
\end{thm}
\noindent \textbf{Proof:} Since $M$ satisfies \eqref{eq4.1} and $L_2$ is nowhere zero, then from the condition (IV) of Theorem \ref{main}, 
$$T Q(\widetilde g, \widetilde S) \equiv 0$$
$$\Rightarrow T \equiv 0 \mbox{ or } Q(\widetilde g, \widetilde S) \equiv 0 \ \ \mbox{ [by Lemma \ref{lem1}].}$$ 
Now if $T$ is identically zero, then from Condition (II), $\overline R (L_1 f \widetilde g + L_2 \widetilde S + L_2 \Omega \widetilde g) \equiv 0$.
$$\Rightarrow \overline R \equiv 0 \mbox{ or } L_1 f \widetilde g + L_2 \widetilde S + L_2 \Omega \widetilde g \equiv 0. \ \ \mbox{ [by Lemma \ref{lem1}]}$$
$$\Rightarrow \mbox{base is flat or fiber is Einstein \ \ \ [as $L_2$ is nowhere zero].}$$
Again if $Q(\widetilde g, \widetilde S)\equiv 0$ then from Lemma \ref{lem2}, the fiber $\widetilde M$ is Einstein. This completes the proof.\\
\indent Now from Theorem \ref{main} we can easily get the characterization for a warped product semisymmetric and various pseudosymmetric type manifolds as follows:
\begin{cor}\label{ssn}
Let $M^n = \overline M^p \times_f \widetilde M^{n-p}$ be a warped product manifold. Then $M$ satisfies $R\cdot R = 0$ if and only if the following conditions hold simultaneously:\\
\indent (I) $\overline R\cdot \overline R = 0,$\\
\indent (II) $T_{as}T_{bc}-T_{ac}T_{bs}+ f T^t_{s}\overline R_{abct} = 0,$\\
\indent (III) $T_{as}\left(\widetilde R_{\eta\alpha\beta\gamma}- f \Delta\widetilde G_{\eta\alpha\beta\gamma}\right) +  T^2_{as}\widetilde G_{\eta\alpha\beta\gamma}= 0,$\\
\indent	(IV) $\widetilde R\cdot \widetilde R = f \Delta Q(\widetilde g,\widetilde R)$.
\end{cor}
\begin{cor}
Let $M^n = \overline M^p \times_f \widetilde M^{n-p}$ be a warped product semisymmetric manifold. Then\\
(i) the base $\overline M$ is semisymmetric,\\
(ii) the fiber $\widetilde M$ is pseudosymmetric,\\
(iii) the fiber $\widetilde M$ is of constant curvature if $T \neq 0$.\\
(iv) $R\cdot T = 0$.
\end{cor}
\begin{cor}\label{psn} (\cite{DD91a}, \cite{Desz93}, \cite{DGJZ16}) 
Let $M^n = \overline M^p \times_f \widetilde M^{n-p}$ be a warped product manifold. Then $M$ satisfies $R\cdot R = L_R Q(g,R)$ if and only if the following conditions hold simultaneously:\\
\indent (I) $\overline R\cdot \overline R = L_R Q(\overline g,\overline R),$\\
\indent (II) $T^t_s \left(\overline R_{abct}-T_{ac}\overline g_{bt}+T_{bc}\overline g_{at}\right) = 
L_R \left(\overline R_{abcs}-T_{ac}\overline g_{bs}+T_{bc}\overline g_{as}\right),$\\
\indent (III) $T_{as}\left(\widetilde R_{\eta\alpha\beta\gamma}- f \Delta\widetilde G_{\eta\alpha\beta\gamma}\right) + f T^2_{as}\widetilde G_{\eta\alpha\beta\gamma}=
 L_R \left[\overline g_{as}\widetilde R_{\eta\alpha\beta\gamma}+ f (T_{as}-\Delta\overline g_{as})\widetilde G_{\eta\alpha\beta\gamma}\right]$,\\
\indent	(IV) $\widetilde R\cdot \widetilde R = \left(f \Delta + f L_R\right) Q(\widetilde g,\widetilde R)$.
\end{cor}
\begin{cor}
Let $M^n = \overline M^p \times_f \widetilde M^{n-p}$ be a warped product pseudosymmetric manifold $(R\cdot R = L_R Q(g,R))$. Then\\
(i) the base and fiber both are pseudosymmetric,\\
(ii) $R\cdot T = L_R Q(\overline g, T)$.
\end{cor}
\begin{cor}\label{rgpsn}
Let $M^n = \overline M^p \times_f \widetilde M^{n-p}$ be a warped product manifold. Then $M$ satisfies $R\cdot R = L Q(S,R)$ if and only if the following conditions hold simultaneously:\\
\indent (I) $\overline R\cdot \overline R = L Q(\overline S,\overline R) - L (n-p) Q(T,\overline R),$\\
\indent (II) $f \widetilde g_{\alpha\beta} \left(T_{as}T_{bc}-T_{ac}T_{bs}+T^t_{s}\overline R_{abct}\right) = 
L \overline R_{abcs}\left(\widetilde S_{\alpha\beta} + \Omega \widetilde g_{\alpha\beta}\right)\\
\indent\hspace{2.6in} +L f \widetilde g_{\alpha\beta}  \left[T_{bc} \left(\overline S_{as}-(n-p)T_{as}\right)-T_{ac} \left(\overline S_{bs}-(n-p)T_{bs}\right)\right]$,\\
\indent (III) $T_{as}\left(\widetilde R_{\eta\alpha\beta\gamma}- f \Delta\widetilde G_{\eta\alpha\beta\gamma}\right) + f T^2_{as}\widetilde G_{\eta\alpha\beta\gamma}=
 L (\overline S_{as}-(n-p) T_{as})(\widetilde R_{\eta\alpha\beta\gamma}- f\Delta \widetilde G_{\eta\alpha\beta\gamma})\\
\indent\hspace{3.3in} + L T_{as} \left[\widetilde g_{\alpha\beta}(\widetilde S_{\gamma\eta} + \Omega \widetilde g_{\gamma\eta}) -\widetilde g_{\alpha\gamma}(\widetilde S_{\beta\eta} + \Omega \widetilde g_{\beta\eta})\right]$,\\
\indent	(IV) $L T_{ab} Q(\widetilde g, \widetilde S)_{\alpha\beta\gamma\delta} =0$ and\\
\indent	(V) $\widetilde R\cdot \widetilde R = \left(f \Delta + \Omega L\right) Q(\widetilde g,\widetilde R) + L Q(\widetilde S,\widetilde R) - L \Delta Q(\widetilde S,\widetilde G)$.
\end{cor}
\begin{cor}
If a warped product manifold $M^n = \overline M^p \times_f \widetilde M^{n-p}$ satisfies $R\cdot R = L Q(S,R)$, then 
$$M= \{x\in M : \overline R|_{\pi(x)} =0\} \cup \{x\in M : T|_{\pi(x)} =0\} \cup \{x\in M : (\widetilde S-\mbox{$\frac{\widetilde \kappa}{n-p}$}\widetilde g)|_{\sigma(x)} =0\}.$$
\end{cor}
\begin{cor}\label{srgpsn} \cite{DD91a}
Let $M^n = \overline M^p \times_f \widetilde M^{n-p}$ be a warped product manifold. Then $M$ satisfies the special Ricci generalized pseudosymmetric condition $R\cdot R = Q(S,R)$ if and only if the following conditions hold simultaneously:\\
\indent (I) $\overline R\cdot \overline R = Q(\overline S- (n-p) T,\overline R),$\\
\indent (II) $f \widetilde g_{\alpha\beta} \left(T_{as}T_{bc}-T_{ac}T_{bs}+T^t_{s}\overline R_{abct}\right) = 
\overline R_{abcs}\left(\widetilde S_{\alpha\beta} + \Omega \widetilde g_{\alpha\beta}\right)\\
\indent\hspace{2.6in} + f \widetilde g_{\alpha\beta}  \left[T_{bc} \left(\overline S_{as}-(n-p)T_{as}\right)-T_{ac} \left(\overline S_{bs}-(n-p)T_{bs}\right)\right]$,\\
\indent (III) $T_{as}\left(\widetilde R_{\eta\alpha\beta\gamma}- f \Delta\widetilde G_{\eta\alpha\beta\gamma}\right) + f T^2_{as}\widetilde G_{\eta\alpha\beta\gamma}=
(\overline S_{as}-(n-p) T_{as})(\widetilde R_{\eta\alpha\beta\gamma}- f\Delta \widetilde G_{\eta\alpha\beta\gamma})\\
\indent\hspace{3.3in} + T_{as} \left[\widetilde g_{\alpha\beta}(\widetilde S_{\gamma\eta} + \Omega \widetilde g_{\gamma\eta}) -\widetilde g_{\alpha\gamma}(\widetilde S_{\beta\eta} + \Omega \widetilde g_{\beta\eta})\right]$,\\
\indent	(IV) $T_{ab} Q(\widetilde g, \widetilde S)_{\alpha\beta\gamma\delta} =0$ and\\
\indent	(V) $\widetilde R\cdot \widetilde R = \left(f \Delta + \Omega \right) Q(\widetilde g,\widetilde R) + Q(\widetilde S,\widetilde R) - \Delta Q(\widetilde S,\widetilde G)$.
\end{cor}
\begin{cor}
If a warped product manifold $M^n = \overline M^p \times_f \widetilde M^{n-p}$ satisfies $R\cdot R = Q(S,R)$, then 
$$M= \{x\in M : \overline R|_{\pi(x)} =0\} \cup \{x\in M : (\widetilde S-\mbox{$\frac{\widetilde \kappa}{n-p}$}\widetilde g)|_{\sigma(x)} =0\}.$$
\end{cor}
\begin{thm}
If a warped product manifold $M^n = \overline M^p \times_f \widetilde M^{n-p}$ satisfies $R\cdot R = Q(S,R)$, then\\
either (i)  the base $\overline M$ is flat\\
or (ii) the fiber $\widetilde M$ is Einstein.
\end{thm}
\noindent \textbf{Proof:} Since the condition $R\cdot R = Q(S, R)$ is a special case of $R\cdot R = L_1 Q(g, R) + L_2 Q(S, R)$, for $L_1 = 0$ and $L_2 =1$. Therefore as $L_2$ is nowhere zero, so from Theorem \ref{thm5.3}, we get our assertion.
\begin{cor}\label{wssn}
Let $M^n = \overline M^p \times_f \widetilde M^{n-p}$ be a warped product manifold. Then $M$ satisfies $W\cdot R = 0$ if and only if the following conditions hold simultaneously:\\
\indent (I) $\overline R\cdot \overline R = \frac{\kappa}{n(n-1)} Q(\overline g,\overline R),$\\
\indent (II) $T^t_s \left(\overline R_{abct}-T_{ac}\overline g_{bt}+T_{bc}\overline g_{at}\right) = 
\frac{\kappa}{n(n-1)} \left(\overline R_{abcs}-T_{ac}\overline g_{bs}+T_{bc}\overline g_{as}\right),$\\
\indent (III) $T_{as}\left(\widetilde R_{\eta\alpha\beta\gamma}- f \Delta\widetilde G_{\eta\alpha\beta\gamma}\right) + f T^2_{as}\widetilde G_{\eta\alpha\beta\gamma}=
\frac{\kappa}{n(n-1)} \left[\overline g_{as}\widetilde R_{\eta\alpha\beta\gamma}+ f (T_{as}-\Delta\overline g_{as})\widetilde G_{\eta\alpha\beta\gamma}\right]$,\\
\indent	(IV) $\widetilde R\cdot \widetilde R = f \left(\Delta + \frac{\kappa}{n(n-1)}\right) Q(\widetilde g,\widetilde R)$.
\end{cor}
\begin{cor}
Let $M^n = \overline M^p \times_f \widetilde M^{n-p}$ be a warped product manifold. Then $M$ satisfies $W\cdot R = L_2 Q(S, R)$ if and only if the following conditions hold simultaneously:\\
\indent $(I)$ $\overline R\cdot \overline R = \frac{\kappa}{n(n-1)} Q(\overline g,\overline R) + L_2 Q(\overline S,\overline R) - L_2 (n-p) Q(T,\overline R),$\\
\indent $(II)$ $f \widetilde g_{\alpha\beta} \left(T_{as}T_{bc}-T_{ac}T_{bs}+T^t_{s}\overline R_{abct}\right) = 
\frac{\kappa}{n(n-1)} f  \widetilde g_{\alpha\beta} \left(\overline R_{abcs}-T_{ac}\overline g_{bs}+T_{bc}\overline g_{as}\right)
\\
\indent\hspace{2.6in} +L_2 \overline R_{abcs}\left(\widetilde S_{\alpha\beta} + \Omega \widetilde g_{\alpha\beta}\right)\\
\indent\hspace{2.6in} +L_2 f \widetilde g_{\alpha\beta}  \left[T_{bc} \left(\overline S_{as}-(n-p)T_{as}\right)-T_{ac} \left(\overline S_{bs}-(n-p)T_{bs}\right)\right]$,\\
\indent $(III)$ $T_{as}\left(\widetilde R_{\eta\alpha\beta\gamma}- f \Delta\widetilde G_{\eta\alpha\beta\gamma}\right) + f T^2_{as}\widetilde G_{\eta\alpha\beta\gamma}=
 \frac{\kappa}{n(n-1)} \left[\overline g_{as}\widetilde R_{\eta\alpha\beta\gamma}+ f (T_{as}-\Delta\overline g_{as})\widetilde G_{\eta\alpha\beta\gamma}\right]\\
\indent\hspace{3.3in} + L_2 (\overline S_{as}-(n-p) T_{as})(\widetilde R_{\eta\alpha\beta\gamma}- f\Delta \widetilde G_{\eta\alpha\beta\gamma})\\
\indent\hspace{3.3in} + L_2 T_{as} \left[\widetilde g_{\alpha\beta}(\widetilde S_{\gamma\eta} + \Omega \widetilde g_{\gamma\eta}) -\widetilde g_{\alpha\gamma}(\widetilde S_{\beta\eta} + \Omega \widetilde g_{\beta\eta})\right]$,\\
\indent	$(IV)$ $L_2 T_{ab} Q(\widetilde g, \widetilde S)_{\alpha\beta\gamma\delta} =0$ and\\
\indent	$(V)$ $\widetilde R\cdot \widetilde R = f \left(\Delta + \frac{\kappa}{n(n-1)} + \Omega L_2\right) Q(\widetilde g,\widetilde R) + L_2 Q(\widetilde S,\widetilde R) - L_2 \Delta Q(\widetilde S,\widetilde G)$.
\end{cor}
\begin{cor}
Let $M^n = \overline M^p \times_f \widetilde M^{n-p}$ be a warped product manifold. Then $M$ satisfies $P\cdot R = 0$ if and only if the following conditions hold simultaneously:\\
\indent (I) $\overline R\cdot \overline R = \frac{1}{n-2} Q(\overline S-(n-p)T,\overline R),$\\
\indent (II) $f \widetilde g_{\alpha\beta} \left(T_{as}T_{bc}-T_{ac}T_{bs}+T^t_{s}\overline R_{abct}\right) = 
\frac{1}{n-2} \overline R_{abcs}\left(\widetilde S_{\alpha\beta} + \Omega \widetilde g_{\alpha\beta}\right)\\
\indent\hspace{2.6in} +\frac{1}{n-2} f \widetilde g_{\alpha\beta}  \left[T_{bc} \left(\overline S_{as}-(n-p)T_{as}\right)-T_{ac} \left(\overline S_{bs}-(n-p)T_{bs}\right)\right]$,\\
\indent (III) $T_{as}\left(\widetilde R_{\eta\alpha\beta\gamma}- f \Delta\widetilde G_{\eta\alpha\beta\gamma}\right) + f T^2_{as}\widetilde G_{\eta\alpha\beta\gamma}=
\frac{1}{n-2} (\overline S_{as}-(n-p) T_{as})(\widetilde R_{\eta\alpha\beta\gamma}- f\Delta \widetilde G_{\eta\alpha\beta\gamma})\\
\indent\hspace{3.3in} + \frac{1}{n-2} T_{as} \left[\widetilde g_{\alpha\beta}(\widetilde S_{\gamma\eta} + \Omega \widetilde g_{\gamma\eta}) -\widetilde g_{\alpha\gamma}(\widetilde S_{\beta\eta} + \Omega \widetilde g_{\beta\eta})\right]$,\\
\indent	(IV) $T_{ab} Q(\widetilde g, \widetilde S)_{\alpha\beta\gamma\delta} =0$ and\\
\indent	(V) $\widetilde R\cdot \widetilde R = f \left(\Delta + \frac{\Omega}{n-2}\right) Q(\widetilde g,\widetilde R) + \frac{1}{n-2} Q(\widetilde S,\widetilde R) - \frac{1}{n-2} \Delta Q(\widetilde S,\widetilde G)$.
\end{cor}
\begin{thm}
If a warped product manifold $M^n = \overline M^p \times_f \widetilde M^{n-p}$ satisfies $P\cdot R = 0$, then
either (i)  the base $\overline M$ is flat 
or (ii) the fiber $\widetilde M$ is Einstein.
\end{thm}
\begin{cor}\label{ppsn}
Let $M^n = \overline M^p \times_f \widetilde M^{n-p}$ be a warped product manifold. Then $M$ satisfies the pseudosymmetric type condition
\be
P\cdot R = L_1 Q(g,R)
\ee
if and only if the following conditions hold simultaneously:\\
\indent (I) $\overline R\cdot \overline R = L_1 Q(\overline g,\overline R) + \frac{1}{n-2} Q(\overline S,\overline R) - \frac{n-p}{n-2}  Q(T,\overline R),$\\
\indent (II) $f \widetilde g_{\alpha\beta} \left(T_{as}T_{bc}-T_{ac}T_{bs}+T^t_{s}\overline R_{abct}\right) = 
L_1 f  \widetilde g_{\alpha\beta} \left(\overline R_{abcs}-T_{ac}\overline g_{bs}+T_{bc}\overline g_{as}\right)
+\frac{1}{n-2} \overline R_{abcs}\left(\widetilde S_{\alpha\beta} + \Omega \widetilde g_{\alpha\beta}\right)\\
\indent\hspace{2.6in} + \frac{f}{n-2} \widetilde g_{\alpha\beta}  \left[T_{bc} \left(\overline S_{as}-(n-p)T_{as}\right)-T_{ac} \left(\overline S_{bs}-(n-p)T_{bs}\right)\right]$,\\
\indent (III) $T_{as}\left(\widetilde R_{\eta\alpha\beta\gamma}- f \Delta\widetilde G_{\eta\alpha\beta\gamma}\right) + f T^2_{as}\widetilde G_{\eta\alpha\beta\gamma}=
 L_1 \left[\overline g_{as}\widetilde R_{\eta\alpha\beta\gamma}+ f (T_{as}-\Delta\overline g_{as})\widetilde G_{\eta\alpha\beta\gamma}\right]\\
\indent\hspace{3.3in} + \frac{1}{n-2} (\overline S_{as}-(n-p) T_{as})(\widetilde R_{\eta\alpha\beta\gamma}- f\Delta \widetilde G_{\eta\alpha\beta\gamma})\\
\indent\hspace{3.3in} + \frac{1}{n-2} T_{as} \left[\widetilde g_{\alpha\beta}(\widetilde S_{\gamma\eta} + \Omega \widetilde g_{\gamma\eta}) -\widetilde g_{\alpha\gamma}(\widetilde S_{\beta\eta} + \Omega \widetilde g_{\beta\eta})\right]$,\\
\indent	(IV) $T_{ab} Q(\widetilde g, \widetilde S)_{\alpha\beta\gamma\delta} =0$ and\\
\indent	(V) $\widetilde R\cdot \widetilde R = \left(f \Delta + f L_1 + \frac{\Omega}{n-2}\right) Q(\widetilde g,\widetilde R) + \frac{1}{n-2} Q(\widetilde S,\widetilde R) - \frac{1}{n-2} \Delta Q(\widetilde S,\widetilde G)$.
\end{cor}
\begin{thm}
If a warped product manifold $M^n = \overline M^p \times_f \widetilde M^{n-p}$ satisfies $P\cdot R = L_1 Q(g,R)$, then either (i)  the base $\overline M$ is flat 
or (ii) the fiber $\widetilde M$ is Einstein.
\end{thm}
\section{\bf Examples}\label{exam}
\textbf{Example 1:} Consider the warped product $M = \overline M\times_f \widetilde M$, where $\overline M$ is an open interval of $\mathbb R$ with the metric $d\overline s^2 = \frac{1}{1 + a (1 + x^1)^2}(dx^1)^2$ in local coordinate $x^1$, $\widetilde M$ is a 4-dimensional manifold equipped with a semi-Riemannian metric
$$d\widetilde s^2 = -(dx^2)^2 + e^{2 x^2} (x^5)^2 (dx^3)^2 + 2e^{2 x^2} dx^3 dx^4 + e^{2 x^2} (dx^5)^2$$
in local coordinates $(x^2, x^3, x^4, x^5)$ and the warping function $f = (x^1 + 1)^2$. We can easily evaluate the local components of necessary tensors of $\widetilde M$. The non-zero components of the Riemann-Christoffel curvature tensor $\widetilde R$ and the Ricci tensor $\widetilde S$ of $\widetilde M$ upto symmetry are
$$\widetilde R_{1213} = \widetilde R_{1414} = -e^{2 x^2}, \ \ -\widetilde R_{2323} = \widetilde R_{2434}=e^{4 x^2},$$
$$\widetilde R_{1212} = -e^{2 x^2} (x^5)^2, \ \ \widetilde R_{2424}=e^{2 x^2} (e^{x^2} x^5-1) (e^{x^2} x^5+1)$$
and
$$\widetilde S_{11}=3, \ \ \widetilde S_{22} = -3 e^{2 x^2} (x^5)^2-1, \ \ \widetilde S_{23}= \widetilde S_{44} = -3 e^{2 x^2}.$$
Scalar curvature of $\widetilde M$ is $(-12)$. Again the non-zero components of $\widetilde R\cdot \widetilde R$, $Q(\widetilde g,\widetilde R)$ and $Q(\widetilde S,\widetilde R)$ are
$$\widetilde R\cdot \widetilde R_{122414}= -\widetilde R\cdot \widetilde R_{142412}=e^{2 x^2}, \ \ 2\widetilde R\cdot \widetilde R_{232424}= -\widetilde R\cdot \widetilde R_{242423}=2 e^{4 x^2},$$
$$-Q(\widetilde g,\widetilde R)_{122414} = Q(\widetilde g,\widetilde R)_{142412}=e^{2 x^2}, \ \ -2Q(\widetilde g,\widetilde R)_{232424}= Q(\widetilde g,\widetilde R)_{242423}=2 e^{4 x^2},$$
$$-Q(\widetilde S,\widetilde R)_{121223}= -2Q(\widetilde S,\widetilde R)_{121424}= 2Q(\widetilde S,\widetilde R)_{122312}= \frac{2}{3}Q(\widetilde S,\widetilde R)_{122414}= -Q(\widetilde S,\widetilde R)_{142412}=2 e^{2 x^2},$$
$$2Q(\widetilde S,\widetilde R)_{232424} = -Q(\widetilde S,\widetilde R)_{242423}=4 e^{4 x^2}.$$
Then we can easily check that $\widetilde M$ satisfies $\widetilde R\cdot \widetilde R = - Q(\widetilde g,\widetilde R)$, i.e.,  $\widetilde M$ is a pseudosymmetric manifold of constant type.\\
\indent Now by a straightforward calculation we can evaluate the components of various necessary tensors corresponding to $M$. The non-zero local components of the Riemann-Christoffel curvature tensor $R$ and the Ricci tensor $S$ of $M$ upto symmetry are
$$R_{1212}=\frac{a (x^1+1)^2}{a (x^1)^2+2 a x^1+a+1}, \ \ R_{1313} = -\frac{a e^{2 x^2} (x^1+1)^2 (x^5)^2}{a (x^1)^2+2 a x^1+a+1},$$
$$R_{1314} = R_{1515} = -\frac{a e^{2 x^2} (x^1+1)^2}{a (x^1)^2+2 a x^1+a+1}, \ \ R_{2323}=a e^{2 x^2} (x^1+1)^4 (x^5)^2,$$
$$R_{2324}= R_{2525}=a e^{2 x^2} (x^1+1)^4, \ \ R_{3434}= -R_{3545}=a e^{4 x^2} (x^1+1)^4,$$
$$-R_{3535}=e^{2 x^2} (x^1+1)^2 (a e^{2 x^2} (x^5)^2 + a e^{2 x^2} (x^1)^2 (x^5)^2+2 a e^{2 x^2} x^1 (x^5)^2+1)$$
and
$$S_{11}=\frac{4 a}{a (x^1)^2+2 a x^1+a+1}, \ \ S_{22}= -4 a (x^1+1)^2,$$
$$S_{33}=4 a e^{2 x^2} (x^5)^2+4 a e^{2 x^2} (x^1)^2 (x^5)^2+8 a e^{2 x^2} x^1 (x^5)^2+1, \ \ S_{34}= S_{55}=4 a e^{2 x^2} (x^1+1)^2.$$
The scalar curvature of $M$ is $20 a$. Now the non-zero components (upto symmetry) of $R\cdot R$, $Q(g,R)$ and $Q(S,R)$ are
$$R\cdot R_{133515} = -R\cdot R_{153513}=\frac{a e^{2 x^2} (x^1+1)^2}{a (x^1)^2+2 a x^1+a+1},$$
$$-R\cdot R_{233525}= R\cdot R_{253523}=a e^{2 x^2} (x^1+1)^4,$$
$$-2R\cdot R_{343535}= R\cdot R_{353534}=2 a e^{4 x^2} (x^1+1)^4,$$
$$Q(g,R)_{133515}= -Q(g,R)_{153513}=\frac{e^{2 x^2} (x^1+1)^2}{a (x^1)^2+2 a x^1+a+1},$$
$$-Q(g,R)_{233525}= Q(g,R)_{253523}=e^{2 x^2} (x^1+1)^4,$$
$$-2Q(g,R)_{343535}= Q(g,R)_{353534}=2 e^{4 x^2} (x^1+1)^4$$
and
$$-Q(S,R)_{121323}= Q(S,R)_{122313}=\frac{a (x^1+1)^2}{a (x^1)^2+2 a x^1+a+1},$$
$$-\frac{1}{2}Q(S,R)_{131334}= -Q(S,R)_{131535}= Q(S,R)_{133413}= \frac{1}{4}Q(S,R)_{133515}$$ $$= -\frac{1}{3}Q(S,R)_{153513}=\frac{a e^{2 x^2} (x^1+1)^2}{a (x^1)^2+2 a x^1+a+1},$$
$$\frac{1}{2}Q(S,R)_{232334}= Q(S,R)_{232535}= -Q(S,R)_{233423}= -\frac{1}{4}Q(S,R)_{233525}$$ $$= \frac{1}{3}Q(S,R)_{253523}=a e^{2 x^2} (x^1+1)^4,$$
$$-2Q(S,R)_{343535} = Q(S,R)_{353534}=6 a e^{4 x^2} (x^1+1)^4.$$
Then from the values of $R\cdot R$ and $Q(g,R)$, we see that $M$ satisfies $R\cdot R = a Q(g,R)$. Therefore $M$ is a warped product pseudosymmetric manifold of constant type. In particular if $a=0$, then $M$ is a warped product semisymmetric manifold. Now $R\cdot R = a Q(g, R) \Rightarrow (R-a G)\cdot R = 0$ and $W=R-\frac{\kappa}{5\times 4} G = R- a G$ hence $W\cdot R = 0$. Therefore $M$ is semisymmetric type due to concircular curvature tensor.
\begin{rem}
We know that fibers of a warped product manifold are totally umbilical submanifold of it. In the above example, $\widetilde M$ is a totally umbilical hypersurface of $M$ for $a\ne 0$ and for $a=0$, the above example is an example of pseudosymmetric totally umbilical hypersurface of a semisymmetric manifold.
\end{rem}
\textbf{Example 2:} \cite{SKppsn} Let $M_2 = \mathbb R \times_f \mathbb R^3$ be the $4$-dimensional connected semi-Riemannian warped product manifold, where the base metric is given by $d\overline s^2 = (1+2e^{x^1})(dx^1)^2$, the fiber metric is the usual Euclidean metric and the warping function $f = (1+2e^{x^1})$. Hence the metric of $M_2$ is given by
\be\label{met2}
ds^2 = (1+2e^{x^1})(dx^1)^2 + (1+2e^{x^1})\left[(dx^2)^2+(dx^3)^2+(dx^4)^2\right].
\ee
Therefore the non-zero components (upto symmetry) of $R$, $S$, $\kappa$ and $P$ are given by
$$R_{1212}= R_{1313}= R_{1414}= -\frac{e^{x^1}}{2 e^{x^1}+1}, \ \ R_{2323}= R_{2424}= R_{3434}= -\frac{e^{2 x^1}}{2 e^{x^1}+1};$$
$$S_{11}=\frac{3 e^{x^1}}{\left(2 e^{x^1}+1\right){}^2}, \ \ S_{22}= S_{33}= S_{44}=\frac{e^{x^1}}{2 e^{x^1}+1}; \ \ \kappa = \frac{6 e^{x^1} (1+e^{x^1})}{(1+2 e^{x^1})^3};$$
$$\frac{1}{2}P_{1221}= \frac{1}{2}P_{1331}= \frac{1}{2}P_{1441}= P_{2323}= -P_{2332}= P_{2424}= -P_{2442} = P_{3434}= -P_{3443}=-\frac{e^{2x^1}-e^{x^1}}{6 e^{x^1}+3}.$$
Using above we can easily calculate the non-zero components (upto symmetry) of $R\cdot R$, $Q(g,R)$, $Q(S,R)$ and $P\cdot R$ as follows:
$$R\cdot R_{122313}= R\cdot R_{122414}= -R\cdot R_{132312}= R\cdot R_{133414}= -R\cdot R_{142412}=-R\cdot R_{143413}=\frac{e^{2 x^1} \left(e^{x^1}-1\right)}{\left(2 e^{x^1}+1\right)^3};$$
$$Q(g,R)_{122313}= Q(g,R)_{122414}= -Q(g,R)_{132312}= Q(g,R)_{133414}= -Q(g,R)_{142412}= -Q(g,R)_{143413}$$$$=e^{x^1} \left(e^{x^1}-1\right);$$
$$Q(S,R)_{122313}= Q(S,R)_{122414}= -Q(S,R)_{132312}= Q(S,R)_{133414}= -Q(S,R)_{142412}= -Q(S,R)_{143413}$$$$=\frac{e^{2 x^1} \left(e^{x^1}-1\right)}{\left(2 e^{x^1}+1\right)^3};$$
$$P\cdot R_{122313}= P\cdot R_{122414}= -P\cdot R_{132312}= P\cdot R_{133414}= -P\cdot R_{142412}= -P\cdot R_{143413}=\frac{2 e^{2 x^1} \left(e^{x^1}-1\right)}{3 \left(2 e^{x^1}+1\right)^3}.$$
In view of above results we see that $M_2$ satisfies the following pseudosymmetric type conditions:\\
(i) $R\cdot R = \frac{e^{x^1}}{\left(2 e^{x^1}+1\right)^3} Q(g,R) = Q(S,R)$,\\
(ii) $P\cdot R = \frac{2 e^{x^1}}{3 \left(2 e^{x^1}+1\right)^3} Q(g,R)$ and\\
(iii) $R\cdot R = Q(g,R) + \left[1- e^{-x^1} \left(2 e^{x^1}+1\right)^3\right]Q(S,R)$.
%
\section{\bf{Conclusions}}
In this present paper we have found out the necessary and sufficient condition for which a warped product semi-Riemannian manifold $M^n = \overline M^p \times_f \widetilde M^{n-p}$ satisfies the pseudosymmetric type condition $R\cdot R = L_1 Q(g,R) + L_2 Q(S,R)$. It is shown that on each point $x\in M$ of such manifold either $\overline R = 0$ or $T = L_1 g$ or $\widetilde S = \frac{\widetilde\kappa}{(n-p)}\widetilde g$. Moreover if $L_2$ is nowhere zero, then either the base $\overline M$ is flat or the fiber $\widetilde M$ is Einstein. As a special case of the main theorem we get the necessary and sufficient condition for which a warped product semi-Riemannian manifold $M^n = \overline M^p \times_f \widetilde M^{n-p}$ satisfies the pseudosymmetric type condition (i) $R\cdot R =0$, (ii) $R\cdot R = L_1 Q(g,R)$, (iii) $R\cdot R = L_2 Q(S,R)$, (iv) $R\cdot R = Q(S,R)$, (v) $P\cdot R = 0$, (vi) $P\cdot R = L Q(g,R)$.\\
\indent It is proved that the base of a semisymmetric warped product manifold is semisymmetric and fiber is pseudosymmetric, whereas both the base and fiber of a pseudosymmetric warped product manifold are pseudosymmetric. It is also proved that on a warped product pseudosymmetric or semisymmetric manifold due to projective curvature tensor either the base is flat or the fiber is Einstein. By using the fact that the fiber of a semisymmetric warped product manifold is pseudosymmetric, we have established an example which ensures that a pseudosymmetric manifold is a totally umbilical hypersurface of a semisymmetric manifold.


\end{document}